\documentclass[a4paper]{article}

\author{Joaquim Ro\'e
   \thanks {Partially supported by CIRIT 1997FI-00141,
         CAICYT PB95-0274 and
         ``AGE-Algebraic Geometry in Europe" contract no. ERB940557.}
         \\
      \small{Departament d'\`Algebra i Geometria, 
           Universitat de Barcelona,} \\
      \small{Gran Via, 585, 
           E-08007, Barcelona.} \\
      \small{e-mail: jroevell@cerber.mat.ub.es} \\
           }

\usepackage{amstex}
\usepackage{amsthm}
\usepackage{xypic}

\newtheorem{Def}{Definition}[section]
\newtheorem{Lem}[Def]{Lemma}
\newtheorem{Cor}[Def]{Corollary}
\newtheorem{Pro}[Def]{Proposition}
\newtheorem{Teo}[Def]{Theorem}

\newcommand{\BF}{\operatorname{BF}}
\newcommand{\Bl}{\operatorname{Bl}}

\renewcommand{\P}{\mathbb{P}}
\renewcommand{\H}{{\cal H}}
\renewcommand{\O}{{\cal O}}

\begin{document}

\title{On the existence of plane curves
   with imposed multiple points}

\maketitle

\begin{abstract}
We prove that a plane curve of degree $d$ with
$r$ points of multiplicity $m$ must have
\begin{gather*}
 d \geq m\,(r-1) \prod_{i=2}^{r-1} \left(1-{i \over{i^2+r-1}} \right) \\
 d>\left(\sqrt{r-1} -{\pi \over 8} \right)\,m
\end{gather*}
\end{abstract}

\section {Introduction}

In \cite{nag} Nagata showed a counterexample to the fourteenth
problem of Hilbert; in his construction, he proved that, for
$n>3$, a plane curve going with multiplicity at least $m$
through $n^2$ points in general position must have degree strictly
bigger than $nm$. Moreover, he conjectured that this result
should also hold for a non-square number of points, that is,
a curve with multiplicity $m$ at $r\geq 10$ points in
general position must have degree strictly bigger
than $\sqrt{r}\,m$.

This conjecture has been proved only in some particular
cases. In \cite{eva}, Evain proves it for $m$ small enough,
concretely for $r>\left( {{8\,m}\over{4\,m-1}}(m+1)\right)^2$.
In the case of irreducible reduced curves, Xu proved
in \cite{xu} the inequalities $d>\sqrt{r}\, m-{1 \over{2\,\sqrt{r-1}}}$
and $d>\sqrt{r-1}\,m$. As far as we know, the best known bound
for the general case is what follows from Nagata's result,
$d> \left[ \sqrt{r} \right]m$, where $[ \cdot ]$ denotes the integral
part.

In this work we prove the inequalities
\begin{gather*}
 d \geq m\,(r-1) \prod_{i=2}^{r-1} \left(1-{i \over{i^2+r-1}} \right) \\
 d>\left(\sqrt{r-1} -{\pi \over 8} \right)\,m
\end{gather*}
for all $r \geq 10$. This
is better than the known bound for
$r$ in any interval $((n+{\pi \over 8})^2+1,(n+1)^2)$,
$n \in {\Bbb{Z}}$.
Our approach is based on a
specialization of the scheme consisting of $r$ points in
general position with multiplicity $m$ to an appropriate cluster
scheme supported at a single point.

We would like to thank the referee for
his/her very helpful suggestions.

\section{Definitions}

Given an algebraic variety $Z$
over an algebraically closed field $k$, and 
a closed subvariety $Z'$ of $Z$, we will write
$b:\Bl(Z,Z') \longrightarrow Z$ for the blowing-up
of $Z$ with center $Z'$.

Let $p_1 \in S_0= \P^2$, $p_2 \in S_1=\Bl(S_0,\{p_1\})$,
\dots, $p_r \in S_{r-1}=\Bl(S_{r-2},\{p_{r-1}\})$ The
set $\{p_1,p_2,\ldots,p_r\}$ is called a \emph{cluster}
(see \cite{cll}) and the sequence $K=(p_1,p_2,\ldots,p_r)$
is an \emph{ordered cluster}. Here we will be concerned
only with ordered clusters and we will call them simply
clusters. Note that some of the points of a cluster
can be identified to
proper points of $\P^2$, whereas others may lie infinitely
near to preceding points. A \emph{system of multiplicities}
for a cluster $K=(p_1, p_2, \ldots, p_r)$ is a sequence
of integers $(m)=(m_1, m_2, \ldots, m_r)$, and a pair
$(K,m)$ where $K$ is a cluster and $(m)$ a system of multiplicities
is called a \emph{weighted cluster}. We review now
briefly some known results on clusters; for the proofs,
refer to \cite{cma}, \cite{cll}, having in mind the minor
change that we do not require all points in a cluster
to be infinitely near to the first one.

Given a weighted cluster, we have an ideal sheaf and
a zero-dimensional subscheme of $\P^2$ associated to it.
Write $S_K=\Bl(S_{r-1},\{p_r\})$ 
and denote by $\pi_K$ the composition $S_K \rightarrow \P^2$
of the blowing-ups of the points of $K$.
Let $E_i$ be the pullback (total transform) in $S_K$
of the exceptional divisor of blowing up $p_i$.
Then the ideal sheaf
$$ \H_{K,m}=\left(\pi_K\right)_*\O_{S_K}(-m_1 E_1 - m_2 E_2 - \cdots - m_r E_r)$$
defines a zero-dimensional subscheme of $\P^2$,
and the stalks of $\H_{K,m}$ are complete ideals
in the stalks of $\O_{\P^2}$.
Conversely, if $\cal I$ is a coherent sheaf of ideals
on $\P^2$ defining a zero-dimensional scheme whose
stalks are complete ideals then there is
a weighted cluster $(K,m)$ such that ${\cal I}=\H_{K,m}$.
We will call such schemes \emph{cluster schemes}.
Remark that a plane curve contains the cluster
scheme defined by $(K,m)$ if and only if
it goes (virtually, as in \cite{cma},\cite{cll})
through $(K,m)$.
This notion has already been
considered by Greuel-Lossen-Shustin in \cite{gls}
(with the name \emph{generalized singularity scheme})
and also by Harbourne in \cite{har} (with the name
\emph{generalized fat point scheme}).

Given two points $p_i$, $p_j$ in a cluster $K$
with $j>i$, we say that $p_j$ is \emph{proximate} to
$p_i$ if and only if $j=i+1$ and $p_j$ lies on the
exceptional divisor $E \subset S_i$ of blowing up $p_i$,
or $j>i+1$ and $p_j$ lies on the \emph{strict} transform
of $E$.
The \emph{proximity inequality} at $p_i$ is
$$ m_i \geq \sum_{p_j \text{ prox. to } p_i} m_j.$$
A cluster satisfying the proximity inequalities at
all its points is called \emph{consistent}.
It happens that different weighted clusters
$(K_1,m^{(1)})$ and $(K_2,m^{(2)})$ define
the same cluster scheme. In this case
$\H_{K_1,m^{(1)}} = \H_{K_2,m^{(2)}}$ and 
we will say that the two clusters are equivalent.
For example, if $p_2$ is infinitely near $p_1$
then the weighted clusters
\begin{eqnarray*}
& K_1 =(p_1) & \quad m^{(1)}=(1) \\
& K_2 =(p_1,p_2) & \quad m^{(2)}=(0,1) 
\end{eqnarray*}
are equivalent.
However, if we ask that $m^{(i)}>0$ for
all $i$ and $(K,m)$ be consistent,
then the cluster scheme determines 
the weighted cluster, but for the ordering of points.

Given an arbitrary weighted cluster $(K,m)$
there is a procedure called \emph{unloading}
(see \cite[4]{cll}, \cite[IV.II]{enr}, or \cite{cma})
which gives a new system of multiplicities $(m')$
such that $(K,m')$ is consistent and 
equivalent to $(K,m)$. In each step of the
procedure, one \emph{unloads} some amount of multiplicity
on a point $p_i$ whose proximity inequality is not satisfied,
from the ponts proximate to it. This means
that there is an integer $n>0$ such that, increasing the
multiplicity of $p_i$ by $n$ and decreasing the
multiplicity of every point proximate to $p_i$
by $n$, the resulting weighted cluster
is equivalent to $(K,m)$ and satisfies the
proximity inequality at $p_i$. In other words,
if $\tilde E_i \subset S_K$ is the strict transform
of the exceptional divisor of blowing-up $p_i$,
$D=-m_1 E_1 - m_2 E_2 - \cdots - m_r E_r$ and
and $\tilde E_i \cdot D <0$ then one chooses $n$ as
the minimal integer with $\tilde E_i \cdot (D-n \tilde E_i) \geq 0$
and replaces $D$ by $D-n \tilde E_i$.
A finite number of unloading steps lead to the
desired equivalent consistent cluster.

Let $T$ be a variety, which for the moment we will think
of as a fixed base for our constructions. Let
$p:X \rightarrow T$ be a smooth morphism of relative
dimension $n$, and let $i:Y \rightarrow X$ be a
smooth embedding over $p$.

Let us consider the diagonal morphism
$\Delta := Id_Y \times_T i$
which makes the following diagram commutative:
$$
\diagramcompileto{diagonal}
   Y \drrto^i \drto^\Delta \ddrto_{Id_Y} \\
   & Y \times_T X \rto_{\quad p_X} \dto^{p_Y} & X \dto^p \\ 
   & Y \rto_{p \circ i} & T
\enddiagram
$$

The image $\Delta(Y)$ is a closed smooth subvariety
isomorphic to $Y$. Consider the blowing-up
$$ \BF(X,Y,T):=\Bl(Y\times_T X,\Delta(Y))
     \overset{b}\longrightarrow Y\times_T X,$$
and the commutative diagram
$$
\diagramcompileto{familia}
    \BF(X,Y,T) \rrto^{\quad p_X \circ b}
               \dto_{p_Y \circ b} & &X \dto \\
    Y \rrto & & T
\enddiagram
$$

We call $\pi=p_X \circ b$ and $q=p_Y \circ b$.
As $\Delta$ is a smooth embedding over $p$,
it follows that $q$ is smooth, of relative
dimension $n$ (see \cite[19.4]{ega}). We call
$$\BF(X,Y,T) \overset{q}\longrightarrow Y$$
the family of blowing up $X$ at the points of $Y$.
We are going to see that the morphism
$\BF(X,Y,T) \overset{\pi}\longrightarrow X$,
makes the fibers of $q$ into ordinary blowing-ups,
hence the name.
Given $y \in Y$, with $p(y)=t$, call
$\BF(X,Y,T)_y=\BF(X,Y,T) \times_Y \{y\}$ and
$X_t=X \times_T \{t\}$. Note that $y \in X_t$.

\begin{Pro}
\label{kleiman}
For every point $y \in Y$, and $t=p(y) \in T$
consider the blowing-up $b: \Bl(X_t,\{y\}) \rightarrow X_t$.
Then there is a unique isomorphism
$$\Bl(X_t,\{y\}) \overset{\psi}\longrightarrow \BF(X,Y,T)_y$$
satisfying $b=\pi|_{\BF(X,Y,T)_y} \circ \psi$.
\end{Pro}
\begin{proof} Follows from \cite[2.4]{klit}, as $\Delta(Y)$ is
obviously a local complete intersection, flat over $Y$.
\end{proof}

\section{Varieties of clusters}
\label{construccio}

Take now $X_{-1}=Spec \, k$, $X_0=\P^2_k$,
$p_0:\P^2_k \rightarrow Spec \, k$, and
define recursively $X_i$, $p_i$ as the blowing-up family
$$
   X_i=\BF(X_{i-1},X_{i-1},X_{i-2})
      \overset{p_i}\longrightarrow X_{i-1}.
$$
The morphisms $p_i$ are
in this case projective and smooth 
of relative dimension 2, so their fibers
are projective smooth surfaces. We have also morphisms
$\pi_i:X_i \rightarrow X_{i-1}$ whose restrictions
to the fibers of $p_i$ are by proposition \ref{kleiman}
the blowing-ups of the points of the fibers of $p_{i-1}$.
To simplify notations, let us say
$\pi_{r,i}=\pi_{i+1} \circ \pi_{i+2} \circ \cdots \circ \pi_r$,
$p_{r,i}=p_i \circ p_{i+1} \circ \cdots \circ p_{r}$.
If there is no confusion possible on $r$, we will also write
$p_i$ for $p_{r,i}$, so $p_i(x)$ is a point in $X_{i-1}$, defined
for all $x$ in $X_r$, $r \geq i$.
For any point $x \in X_i$, we will call
$S_x=(X_i)_{p_i(x)} = X_i \times_{X_{i-1}} \{p_i(x)\}$
the surface containing $x$.
Recall that for any cluster $K$, 
$\pi_K:S_K \rightarrow \P^2$ is the composition
of the blowing-ups of the points in $K$.

The following proposition makes the set of all
clusters with $r$ points into an
algebraic variety.

\begin{Pro}
\label{clustersfam}
For every $r \geq 1$ there is a bijection
$$ X_{r-1} \overset{K}\longrightarrow
   \left\{ \text{\emph{clusters of }}
            r \text{ \emph{points}} \right\}
$$
and, for every $x\in X_{r-1}$, a unique isomorphism
$\psi_x:S_{K(x)}\rightarrow (X_r)_x$ such that
$\pi_K=\pi_{r,0}|_{(X_r)_x} \circ \psi_x$.
\end{Pro}
\begin{proof}
Follows from \cite[1.2]{hasun}, since there is
an obvious bijection
\begin{align*}
 \left\{ \text{\emph{ordered blowing-ups at }}
            r \text{ \emph{points}} \right\}
                    &  \longrightarrow
 \left\{ \text{\emph{ordered clusters of }}
            r \text{ \emph{points}} \right\} \\
 S_K &\longmapsto K
\end{align*}
\end{proof}

Notice that the ordering of points in clusters
is essential in proposition \ref{clustersfam}.
If two clusters differing only in the order
of points were considered equal,
as in \cite{cll}, then
injectivity would fail. 
From now on identify the set of clusters of $r$ 
points to the variety $X_{r-1}$.

For every pair of integers $1 \leq i < j \leq r $
there is a subset of $X_{r-1}$ containing exactly those
clusters $K=(x_1, x_2, \ldots, x_r)$ for which
$x_j$ is proximate to $x_i$. It can be proved
that these subsets are constructible subsets of $X_{r-1}$;
we will focus on some of them which are
irreducible closed varieties.

Call $F_i$ the exceptional divisor of
$$ X_i \overset{b_i} {\longrightarrow} X_{i-1} \times_{X_{i-2}} X_{i-1}.$$
Because of proposition \ref{kleiman}
the pullback of $F_i$ to $(X_i)_{p_i}$
is the exceptional divisor $E_i$ of blowing
up $p_i$ in $S_{p_i}$.
It is clear that $p_i(K)$
is proximate to $p_{i-1}(K)$ if and
only if $p_i(K)\in F_{i-1}$.
So there is a closed subvariety
$$ Y_{r-1}:=\bigcap_{i=2}^r p_i^{-1}(F_{i-1}) \subset X_{r-1} $$
containing exactly those clusters $K$
for which $p_{i+1}(K)$ is proximate to $p_{i}(K)$
for all $i$. It is also clear that
$p_{r-1}(Y_{r-1})=Y_{r-2}$, if we allow $Y_0=\P^2$.

\begin{Lem}
\label{tirallonga}
For all $r$, there is a closed immersion
$$ \BF(X_{r-1},Y_{r-1},X_{r-2})
   \overset{i}\longrightarrow X_r$$
such that $Y_r$ is the image of the exceptional
divisor $F'_r$ of 
$$ \BF(X_{r-1},Y_{r-1},X_{r-2})
   \overset{b}\longrightarrow
   X_{r-1} \times_{X_{r-2}} Y_{r-1}.$$
\end{Lem}

\begin{proof}
The closed immersion $i$ is the strict transform
of the closed immersion
$$ Y_{r-1} \times_{X_{r-2}} X_{r-1} 
   \longrightarrow X_{r-1} \times_{X_{r-2}} X_{r-1}
$$
(see \cite[II,7.15]{hag}).
By definition of the $Y_r$ we know that
$Y_r=F_r \cap p_r^{-1}(Y_{r-1}),$
and obviously $i(F'_r) \subset F_r$, and
$$ (p_r \circ i)(\BF(X_{r-1},Y_{r-1},X_{r-2}))=Y_{r-1},$$
so $i(F'_r) \subset Y_r$.
on the other hand, if $y_r\in Y_r$
then $p_r(y_r)=y_{r-1}\in Y_{r-1}$, so
$$ y_r \in S_{y_r} \cong \Bl(S_{y_{r-1}},\{y_{r-1}\})
   \cong \BF(X_{r-1},Y_{r-1},X_{r-2})_{y_{r-1}},$$
which implies $y_r\in i(F'_r)$. So $Y_r \subset i(F'_r)$,
and the proof is complete.
\end{proof}

\begin{Cor}
For all $r$, $Y_r$ together with the restricted morphism
$p_r:Y_r \rightarrow Y_{r-1}$ is a $\P^1-$bundle, and $Y_r$
is irreducible.
\end{Cor}

To deal with the proximity relations between points
$p_i$ and $p_j$ where $j>i+1$ we need some control
on the strict transforms of the exceptional divisor
of blowing up $p_i$. In contrast to what we have
seen in the case $j=i+1$, there is no variety
$\tilde F_i \subset X_{j-1}$ whose pullback
to $(X_{j-1})_{p_{j-1}(K)}$ is the desired strict transform
for all $K$. To overcome this difficulty we 
restrict ourselves to clusters in $Y_{r-1}$
and define varieties $D_{i,j} \subset X_{j-1}$
whose pullback to $(X_{j-1})_{p_{j-1}(K)}$ is the strict
transform of the exceptional divisor
of blowing up $p_i(K)$ if $p_{j-1}(K)$ is proximate
to $p_i(K)$ and empty in any other case. Let first
$$D'_{i,i+1}=D_{i,i+1}=Y_i.$$
Suppose now we have defined $D_{i,j-1} \subset X_{j-2}$
and $D'_{i,j-1}=D_{i,j-1}\cap Y_{j-2}$,
such that the morphism $p_{j-2}|_{D_{i,j-1}}$
is smooth of relative dimension 1 (observe that
for $D_{i,i+1}=Y_i$ this is so).
As there is a closed immersion
$D_{i,j-1} \times_{X_{j-3}} D'_{i,j-1} \rightarrow
   X_{j-2} \times_{X_{j-3}} X_{j-2}$
there is also a closed immersion (its
strict transform)
$$ D_{i,j}=\BF \left( {D_{i,j-1}, D'_{i,j-1}}, X_{j-3} \right)
   \overset{i}\longrightarrow X_{j-1}
$$
which we take as the definition of $D_{i,j}$.
Moreover, as $p_{j-2}|_{D_{i,j-1}}$ is
smooth of relative dimension 1, $\Delta(D_{i,j-1})$ has
codimension 1 in $D_{i,j-1} \times_{X_{j-3}} D_{i,j-1}$
and
$$ \BF \left( {D_{i,j-1}, D'_{i,j-1}}, X_{j-3} \right)
   \overset{b}\longrightarrow D_{i,j-1} \times_{X_{j-3}} D'_{i,j-1}$$
is an isomorphism. We have 
$$ D'_{i,j}=D_{i,j} \cap Y_{j-1} \cong \Delta(D'_{i,j-1}) \subset
   \BF \left( {D_{i,j-1}, D'_{i,j-1}}, X_{j-3} \right).$$
So $D'_{i,j}$ is isomorphic
to $D'_{i,j-1}$, and $p_{j-1}|_{D_{i,j}}$ is smooth of
relative dimension 1.

We will call $(i,j)-$proximity variety the subvariety
$P_{i,j} = p_j^{-1} (D'_{i,j}) \subset Y_{r-1}$.

\begin{Lem}
\label{tancats}
In a cluster $K \in Y_{r-1}$ the points $p_{i+1},p_{i+2},\ldots,p_j$
are proximate to $p_i$ if and only if $K$ lies
in the $(i,j)-$proximity variety. Furthermore, the
proximity varieties are irreducible and there are inclusions
$$ P_{i,i+1} \supset P_{i,i+2} \supset \cdots \supset P_{i,r}.$$
\end{Lem}
\begin{proof}
The first part will clearly be proved if we show that
$$D_{i,j} \times_{X_{j-1}} \{p_{j-1}\} \subset S_{p_j}$$
is the strict transform of $E_i$. This comes out easily by induction
on $j-i$. For $j-i=1$, it is immediate by proposition
\ref{kleiman}. For $j-i>1$,
proposition \ref{kleiman} gives
$$D_{i,j} \times_{X_{j-1}} \{p_{j-1}\} =
   \Bl(D_{i,j-1}\times_{X_{j-2}} \{p_{j-2}\},p_{j-1}),$$
that is, the strict transform in $S_{p_j}$
of $D_{i,j-1}\times_{X_{j-2}} \{p_{j-2}\}$, which
by the induction hypothesis is the
strict transform of $E_i$ in $S_{p_{j-1}}$, so we are done.

From their own definition, the $D'_{i,j}$ are all
isomorphic to $Y_i$, which is irreducible.
Induction on $r-j$ gives the irreducibility
of the $P_{i,j}$. Indeed, if
$P^{(r-1)}_{i,j}= (p_j|_{Y_{r-2}})^{-1}(D_{i,j})$
is irreducible then its preimage by 
$p_{r-1}|Y_{r-1}$ must be irreducible also,
because $Y_{r-1} \rightarrow Y_{r-2}$ is a
projective space bundle.

The inclusions between the $P_{i,j}$ are clear,
from the first part of the lemma.
\end{proof}

Lemma \ref{tancats} shows that there are subsets
$U_{1,i}$ open and dense in $P_{1,i}$ which contain
all those clusters $K$ with
\begin{itemize}
 \item $p_j(K)$ proximate to $p_{j-1}(K)$, $j=2,\ldots,r$
 \item $p_j(K)$ proximate to $p_1(K)$, $2 \leq j \leq i$
\end{itemize}
and no other proximity relations.

\begin{Lem}
\label{descarrega}
Let $(m)=(m_1,m_2,\ldots,m_r)$ be a system of multiplicities,
and call $M=\sum_{j=2}^r m_j$. 
Define $\alpha_i={{i-1}\over{r-1}}$ and
$\beta_i= 1- {{i-1} \over {(i-1)^2+r-1}}$. 
Suppose that for some $i \in \{2,3,\ldots, r\}$
and $A \in {\Bbb R}$ the inequalities
\begin{gather*}
{{(i-2)\,m_1 + M}\over{(i-2)\,\alpha_{i-1} + 1}} \geq A, \\
m_1 \geq \alpha_{i-1}\,A,
\end{gather*}
are satisfied. Then there is a system of
multiplicities $(m')$ which is equivalent to $(m)$ for
all clusters in $U_{1,i}$ and satisfies
\begin{gather}
\label{desig1}
{{(i-1)\,m'_1 + M'}\over{(i-1)\,\alpha_{i} + 1}} \geq \beta_i \,A,
    \\
\label{desig2}
m'_1 \geq \alpha_{i} \, \beta_i \, A. 
\end{gather}
\end{Lem}

\begin{proof}
We know that for a given cluster of r points $K$
there is a system of multiplicities $(m')$, 
consistent and equivalent to $(m)$, which
is obtained from $(m)$ by the unloading procedure.
The unloading procedure depends only on the multiplicities
and the proximity relations, and so it is the same
for all clusters in $U_{1,i}$.

Due to the proximity relations which hold for the
points of a cluster in $U_{1,i}$, when an unloading
step is applied to the point $p_j$, $1<j<r$
the only point whose multiplicity is decreased
is $p_{j+1}$, so $m_1$ and $M$ remain unchanged.
When an unloading step is applied to $p_1$, the
points whose multiplicity is decreased
are $\{p_2,p_3,\ldots,p_i\}$, so if $m_1$ is
increased by $n$, $M$ is decreased by
$(i-1)\,n$. In both cases, the quantity
$ (i-1)\,m_1 + M $ remains the same. 
When an unloading step is applied to $p_r$,
which happens only when its multiplicity
has become negative, then one replaces it by zero,
so $ (i-1)\,m_1 + M $ might increase, but
does never decrease.
After the complete unloading procedure we get
\begin{gather*}
 (i-1)\,m'_1 + M' \geq (i-1)\,m_1 + M = (i-2)\,m_1 + M + m_1 \geq \\
 \geq ((i-2)\,\alpha_{i-1} + 1) A + \alpha_{i-1}\,A =
((i-1)\,\alpha_{i} + 1) \, \beta_i \,A.
\end{gather*}
This proves (\ref{desig1}). To see (\ref{desig2}), 
we multiply this inequality by $\alpha_i$, so we get
$$
\alpha_i\,((i-1)\,m'_1 + M') \geq 
   ((i-1)\,\alpha_{i} + 1) \, \alpha_i\, \beta_i \,A.
$$
On the other hand, as $(m')$ is consistent, $(K,m)$
must satisfy all the proximity inequalities, and these
imply easily
$$m'_1 - \alpha_i \,M' \geq 0.$$
If we add both inequalities, we obtain (\ref{desig2}).
\end{proof}

\section{The bound}

Let $F_i^{(r)}$ be the pullback of $F_i$ by
$ \pi_{r,i}:X_r \longrightarrow X_i$.
Let $\left[ F_0 \right ]^{(r)}$ be the
pullback to $X_r$ by $\pi_{r,0}$
of the class of a line in $\mathbb{P}^2 $.
For any cluster $K\in X_{r-1}$ and $i>0$,
the pullback to the surface $S_K$ of $F_i^{(r)}$ by
the inclusion is obviously the same as the
pullback $E_i$ of the class of the exceptional divisor 
of blowing up $p_i$ in $S_{p_i(K)}$ by
$\pi_{r,i}|_{S_K}$. Similarly, the pullback
of $\left[ F_0 \right]^{(r)}$ to $S_K$ is
the same as the pullback $\left[ E_0 \right]$
of the class of a line by $\pi_{r,0}|_{S_K}$.
All together, we have
\begin{equation}
\label{restriccio}
 \O_{X_r} (F_i^{(r)}) \otimes_{X_{r-1}} k(K)=
   \O_{S_K} (E_i)
\end{equation}
for all $i$.
Given an integer $d$ we define 
$$ {\cal J}_{d,m}={\cal O}_{X_r}
   (d F_0^{(r)}-m_1 F_1^{(r)}-m_2 F_2^{(r)}-\cdots-m_r F_r^{(r)}).$$
Equality (\ref{restriccio}) and the projection
formula show that,
for every cluster $K\in X_{r-1}$, 
$$\H_{K,m}(d)=\H_{K,m} \otimes \O_{\P^2}(d)=
   \left(\pi_K\right)_*({\cal J}_{d,m}\otimes_{X_{r-1}} k(K))$$ and
$H^0 \left( \H_{K,m}(d) \right) \cong
  H^0 \left( {\cal J}_{d,m}\otimes_{X_{r-1}} k(K)\right).$

In our specialization, we start from a cluster $K$ consisting
of $r$ points in general
position, to specialize it, step by step, to the closed
subvarieties $P_{1,i}$. We obtain
the following theorem:

\begin{Teo}
\label{elteorema}
If a plane curve of degree $d$ passes with multiplicity
$m$ through $r$ points in general position, then
\begin{equation}
\label{condicio}
 d \geq m \,(r-1) \prod_{i = 2}^{r-1} \left(
      1 - {i\over {{i^2} + r-1}} \right)
\end{equation}
\end{Teo}
\begin{proof}
Let $\cal J$ and $\cal H$ be the sheaves defined
above. We start from
the system of multiplicities $(m)=(m, m, \ldots, m)$.
We have to prove that for general $K \in X_{r-1}$, the inequality
$$
H^0(({\cal J}_{d,m})\otimes_{X_{r-1}} k(K)) \,\cong\,
    H^0({\cal H}_{K,m} (d)) \neq 0
$$
implies (\ref{condicio}), so assume this inequality holds
for general $K$.
As $X_r \rightarrow X_{r-1}$ is smooth,
the invertible sheaf ${\cal J}_{d,m}$ is flat
over $X_{r-1}$, so by the semicontinuity theorem
\cite[III,12.8]{hag} we have
$$ H^0({\cal H}_{K,m} (d)) \neq 0$$
for all $K \in P_{1,i}$ and any $i$.

Now for $K \in P_{1,3}$ the system of multiplicities
$(m)$ is not consistent. We can find
by unloading multiplicities a consistent system
$(m^{(3)})$ which is equivalent to $(m)$ for
all clusters in $U_{1,3}$. Applying lemma \ref{descarrega}
with $(m)=(m,m,\ldots,m)$, $M=(r-1)\,m$ and $i=3$, we have
$$
{{(i-2)\,m_1+M} \over {(i-2) \alpha_{i-1}+1}} \,=\,
{{m+(r-1)\,m} \over {\alpha_2+1}} \,=\, m\,(r-1) 
$$                                                 
so we can take $A=A_2=m\,(r-1)$ and the lemma gives
\begin{gather*}
{{2\,m^{(3)}_1 + M^{(3)}} \over {2\,\alpha_{3} + 1}} \geq
   m\,(r-1)\, \beta_3 \\
m^{(3)}_1 \geq m\,(r-1)\, \alpha_{3} \, \beta_3.  
\end{gather*}

As $(m^{(3)})$ is equivalent to $(m)$ for
all clusters in $U_{1,3}$, we have
$$
H^0(\H_{K,m^{(3)}} (d))=H^0(\H_{K,m}(d)) \neq 0
$$
if $K \in U_{1,3}$. As $U_{1,3}$ is open and dense
in $P_{1,3}$, and $P_{1,4}\subset P_{1,3}$, the
semicontinuity theorem applied to the new
sheaf ${\cal J}_{d,m^{(3)}}$ implies
$$ H^0({\cal H}_{K,m^{(3)}} (d)) \neq 0$$
for all $K \in P_{1,4}$.
The new system of multiplicities need not be
(but in fact could be) consistent
for $K \in P_{1,4}$. In any case we can find a new
system $(m^{(4)})$ (which could be equal to
$(m^{(3)})$) to use here.
We apply lemma \ref{descarrega} to the new situation,
with $A_3=m\,(r-1)\,\beta_3$, and we obtain
\begin{gather*}
{{3\,m^{(4)}_1 + M^{(4)}} \over {3\,\alpha_{4} + 1}} \geq
   m\,(r-1)\, \beta_3 \beta_4\\
m^{(4)}_1 \geq m\,(r-1)\, \alpha_{4} \, \beta_3 \, \beta_4. 
\end{gather*}

Iterating the process we finally
get a system $(m^{(r)})=(m_1^{(r)},m_2^{(r)}, \cdots m_r^{(r)})$,
with
$$
m^{(r)}_1 \geq m\,(r-1)\, \alpha_{r} \, \prod_{i = 3}^{r} \beta_i \,=\,
 m \,(r-1) \prod_{i = 2}^{r-1} \left(
      1 - {i\over {{i^2} + r-1}} \right)
$$
and 
$$H^0 ({\cal H}_{K,m^{(r)}}(d)) \neq 0$$
for all $K \in P_{1,r}$. It is clear that
this implies $d \geq m^{(r)}_1$.
\end{proof}

The reader may notice that the
proof of theorem \ref{elteorema} 
is valid for any divisor class on an irreducible
smooth projective surface $S$, except for the
last step, namely $d \geq m^{(r)}_1$, which
assumes $C \subset S = \P^2$.
The specialization of a set of multiple
points to a cluster scheme containing
a point of multiplicity
$ m' \geq m \,(r-1) \prod_{i = 2}^{r-1} \left(
      1 - {i\over {{i^2} + r-1}} \right)$
holds thus on any such surface.

\section{A calculation}
\label{last}

The aim of this section is to compare the
bound of theorem \ref{elteorema} with Nagata's
conjecture (which reads $ d > m \sqrt{r}$), and
with previously known results.
We obtain the following:

\begin{Pro}
\label{producte}
Let $n \geq 9$ be a natural number. Then
$$
   n\, \prod_{i = 2}^{n} \left(1 - {i\over {{i^2} + n}} \right)
   > \sqrt{n}-{\pi \over 8}.
$$
\end{Pro}
This has an immediate corollary:
\begin{Cor}
If a plane curve of degree $d$ passes with multiplicity $m$
through $r \geq 10$ points in general position, then
$$
 d \, > \, m \, \left( \sqrt {r-1} - {\pi \over 8}\right)
$$
\end{Cor}

\begin{proof}[Proof of proposition \ref{producte}]
The goal is to bound
$$b=n\prod_{i=2}^n \left(1-{i \over {i^2+n}} \right) =
    n\prod_{i=1}^{n-1} \left(1-{i \over {i^2+n}} \right) =
     n \prod_{i=1}^{n-1} {{n+i^2-i} \over {i^2+n}}
$$
below. This can be rewritten as
$$
n \ {{(1/n) \prod_{i=2}^{n} (n+i^2-i)} \over
     {\prod_{i=1}^{n-1} (i^2+n)}}=
\prod_{i=1}^{n-1} {{n+(i+1)^2-(i+1)} \over {i^2+n}}=
\prod_{i=1}^{n-1} \left( 1+{{i} \over {i^2+n}} \right) .
$$
We thus have
$$
b^2=n \ \prod_{i=1}^{n-1} \left(1-{i \over {i^2+n}} \right) \
        \prod_{i=1}^{n-1} \left(1+{i \over {i^2+n}} \right)=
    n \ \prod_{i=1}^{n-1}
            \left(1-\left({i \over {i^2+n}} \right)^2 \right).
$$
Let $1-\epsilon=\prod_{i=1}^{n-1}(1-(i/(i^2+n))^2)$ and let
$1+\delta=\prod_{i=1}^{n-1}(1+(i/(i^2+n))^2)$. Then $\epsilon=
1-\prod_{i=1}^{n-1}(1-(i/(i^2+n))^2)$ and
$\delta=-1+\prod_{i=1}^{n-1}(1+(i/(i^2+n))^2)$ both involve the 
same terms, except that they occur with signs in $\epsilon$, so 
$0<\epsilon<\delta$. Thus $1-\delta<1-\epsilon$, and so $b^2>n(1-\delta)$.

We can bound $b^2$ (and hence $b$) below by bounding $1+\delta$ 
(and hence $\delta$) above. But $\log (1+x) \le x$ so 
$\log \prod_{i=1}^{n-1}(1+(i/(i^2+n))^2)\le 
\sum_{i=1}^{n-1}(i/(i^2+n))^2$. 

The Fourier series for $\sinh \sqrt{n}x$ on $[-\pi,\pi]$
is
$${2 \over \pi} (\sinh \sqrt{n}\pi)
\sum_{i\ge1}(-1)^i {-i \over {i^2+n}} \ {\sin ix}
$$
so Parseval's identity gives
$$\left( {\pi \over {2\sinh \sqrt{n}\pi}} \right)^2
{1 \over \pi}
\int_{-\pi}^\pi \sinh^2 \sqrt{n}x \ dx
=\sum_{i\ge1} \left({i \over {i^2+n}} \right)^2 .
$$
The integral can be exactly evaluated; we get 
$$
\int_{-\pi}^\pi \sinh^2 \sqrt{n}x \ dx =
-\pi+{1 \over{2\sqrt{n}}}\sinh 2\sqrt{n}\pi.
$$
Also,
$$
\sum_{i\ge n} \left({i \over {i^2+n}} \right)^2 \ge
\sum_{i\ge n} \left({i \over {i^2+i}} \right)^2 \ge
\int_n^\infty \left({1 \over {x+1}} \right)^2 dx =
{1 \over {n+1}}.
$$
Thus we have
\begin{gather*}
\sum_{i=1}^{n-1} \left( {i \over {i^2+n}} \right)^2 \le
 \left({ \pi \over {2\sinh \sqrt{n}\pi}} \right)^2
 {1 \over \pi} \left(-\pi+{1 \over {2\sqrt{n}}}
 \sinh 2\sqrt{n}\pi \right)-{1\over{n+1}} \le \\
\le {\pi \over {8\sqrt{n}}}
 {{\sinh 2\sqrt{n}\pi} \over {\sinh^2 \sqrt{n}\pi}}
  -{1\over{n+1}}.
\end{gather*}
Define $t=e^{\sqrt{n}\pi}$, so
\begin{gather*}
 {{\sinh 2\sqrt{n}\pi} \over { 2 \sinh^2 \sqrt{n}\pi}} =
 \ {{t+1/t} \over {t-1/t}} = 
 \left(1+{1 \over t^2}\right) \left(1+{1 \over t^2} +
            {1 \over t^4}+\cdots\right) = \\ =
  \left(1+{1 \over t^2}\right) \left(1+{1 \over t^2}
            \left(1+{1 \over t^2} + \cdots\right)\right) \le 
 \left(1+{1 \over {t^2}} \right)
      \left(1+{1.5 \over {t^2}} \right) \le
      \left(1+{3 \over {t^2}} \right)
\end{gather*}
because $n\ge 9$, and
$$
\sum_{i=1}^{n-1} \left( {i \over {i^2+n}} \right)^2 \le
 {\pi \over {4\sqrt{n}}} + {1 \over t^2} - {1\over{n+1}}.
$$
But $e^{\sqrt{n}\pi}\ge 3n$ (look at the tangent line to
$e^{\sqrt{n}\pi}$ at $n=9$), so $e^{2\sqrt{n}\pi}\ge 3n^2\ge (n+2)(n+1)$
hence $1/t^2\le 1/((n+2)(n+1))$;
therefore
$$
\sum_{i=1}^{n-1} \left( {i \over {i^2+n}} \right)^2 \le
{\pi \over {4\sqrt{n}}} - {1\over{n+2}}.
$$
This means $\delta\le -1+e^{\pi/(4\sqrt{n})-1/(n+2)}$,
hence
\begin{gather*}
b^2\ge n(1-\delta) \ge n(2-e^{\pi/(4\sqrt{n})-1/(n+2)})
= \\ =
n \left( 1-{\pi \over {4\sqrt{n}}}+ {1 \over {n+2}} -
 {1  \over {2!}} \left( {\pi \over {4\sqrt{n}}}-
 {1 \over {n+2}} \right)^2  -\cdots \right) \ge \\ \ge
n \left(1-{\pi \over {4\sqrt{n}}}+ {1 \over {n+2}} -
 {1 \over {2!}}\left( {\pi \over {4\sqrt{n}}} \right)^2
  -\cdots \right) 
\end{gather*}
and by comparison with a geometric series, this last is
at least as big as
$$
n-{{\pi\sqrt{n}}\over 4}+{n \over {n+2}}-n{{2u^2}\over {1-u}},
$$
where $u= (\pi/(4\sqrt{n}))/2\le \pi/24$, so
$1/(1-u)\le 1/(1-(\pi/24))\le 1.2$, so 
$-n2u^2/(1-u)\ge -2.4nu^2\ge -.4$; i.e., 
$b^2\ge n-\pi\sqrt{n}/4+n/(n+2)-.4\ge n-\pi\sqrt{n}/4+9/11-.4$.
Finally, $(\sqrt{n}-\pi/8)^2=n-\pi\sqrt{n}/4 +\pi^2/64$,
and $9/11-.4>\pi^2/64$, so $b\ge \sqrt{n}-\pi/8$, as required.
\end{proof}

\end{document}